\documentclass[11pt]{article}

\usepackage{amssymb}
\usepackage{amsmath}
\usepackage{graphicx}

\newtheorem{cor}{Corollary}
\newtheorem{theo}{Theorem}

\newtheorem{lem}{Lemma}
\def \UN{\hbox{ \rm 1\hskip -3.2pt I}}

\def \R{ \mathbb R}

\def \E{\hbox{\rm E}}

\def \P{\hbox{\rm P}}
\def \Var{\hbox{\rm Var}}
\def \Cov{\hbox{\rm Cov}}
\def \Fle { \leadsto}
\def \C  {\mathcal{C}_{t,j}}
\def \ch {\widehat{\mathcal{C}}_{t,j}}

\def \c1 { \mathcal{C}^1 }
\begin{document}

\title{A general expression for the distribution of the maximum of a Gaussian field and
the approximation of the tail\footnote{This work was supported by ECOS
program U03E01.}}

\author{Jean-Marc Aza\"{\i}s \thanks{%
Laboratoire de Statistique et Probabilit\'es. UMR-CNRS C5583
Universit\'e Paul Sabatier. 118, route de Narbonne. 31062 Toulouse Cedex
4. France.}, azais@cict.fr \and Mario Wschebor \thanks{%
Centro de Matem\'atica. Facultad de Ciencias. Universidad de la
Rep\'ublica. Calle Igua 4225. 11400 Montevideo. Uruguay.},
wschebor@cmat.edu.uy}
 \maketitle

\bigskip
AMS subject classification: Primary 60G70 Secondary 60G15

\emph{Short Title: }Distribution of the Maximum.

\emph{Key words and phrases:} Gaussian fields, Rice Formula, Euler-Poincar\'e Characteristic,
Distribution of the Maximum, Density of the Maximum, Random Matrices.

\begin{abstract}

We study the probability distribution $F(u)$ of the maximum
 of smooth Gaussian fields defined on compact subsets of $\R^d$  having some geometric regularity.

 Our main result is a general expression for the density of $F$. Even though this is an implicit formula,
 one can deduce from it explicit bounds for the density, hence for the distribution, as well as improved
 expansions
 for $ 1-F(u)$ for large values of $u$.

 The main tool is the Rice formula for the moments of the number of roots of a random
 system of equations over the reals.

 This method enables also to study second order properties of the expected
 Euler Characteristic approximation using only elementary arguments
 and to extend these kind of results  to some interesting classes
 of Gaussian fields.
 We obtain more precise results for the "direct method" to compute
 the distribution of the maximum, using spectral theory of
GOE random matrices.

\end{abstract}

\section{Introduction and notations}

Let $\mathcal{X}= \{X(t): t\in S \}$ be
 a real-valued random
 field defined on some parameter
set $S$ and $M:= \sup_{t\in S}
X(t)$ its supremum. \medskip

The study of the probability distribution of the random variable
$M$, i.e. the function $F_M(u) := \P \{ M \leq u \}$ is a
classical problem in probability theory. When the process is
Gaussian, general inequalities allow to give bounds on $ 1- F_M(u)
= \P\{ M>u\}$ as well as asymptotic results for $u \to +\infty$.
A partial account of this well established theory, since the
founding paper by Landau and Shepp  \cite{landau} should contain -
among a long list of contributors - the works of Marcus and Shepp
\cite{marcus}, Sudakov and Tsirelson \cite{sudakov}, Borell
\cite{borell75} \cite{borell03}, Fernique \cite{fernique}, Ledoux
and Talagrand \cite{ledouxt}, Berman \cite{berman85} \cite
{berman92}, Adler\cite{adler90}, Talagrand \cite{talagrand} and
Ledoux\cite{ledoux}. \medskip

During the last fifteen years, several methods have been introduced
with the aim of obtaining more precise results than those arising
from the classical theory, at least under certain restrictions on
the process $ \mathcal{X}$, which are interesting from the point of
view of the mathematical theory as well as in many significant
applications. These restrictions include the requirement the domain
$S$ to have certain finite-dimensional geometrical structure and the
paths of the random field to have a certain regularity.

Some examples of these contributions are the double sum method by
Piterbarg \cite{piter96}; the Euler-Poincar\'e Characteristic (EPC)
approximation, Taylor, Takemura  and Adler \cite{TTA}, Adler and
Taylor \cite{AT};  the tube method, Sun \cite{sun} and the well-
known Rice method, revisited by Aza\"{\i}s and Delmas \cite{AD},
Aza\"{\i}s and Wschebor \cite{AWR}.
 See also Rychlik \cite{rychlik} for numerical computations.

 \medskip

The results in the present paper are based upon Theorem \ref{exact}
which is an extension of Theorem 3.1 in Aza\"{\i}s and Wschebor
\cite{AW} allowing to express the density $p_M$ of $F_M$  by means
of a general formula. Even though this is an exact formula, it is
only implicit as an expression for the density, since the relevant
random variable $M$ appears in the right-hand side. However, it can
be usefully employed for various purposes. \medskip

{\bf First}, one can use Theorem \ref{exact} to obtain bounds for
$p_M(u) $ and thus for  $\P\{M>u\}$  for { \bf every} $u$ by means
of replacing some indicator function in (\ref{densm}) by the
condition that the normal derivative is "extended outward" (see
below for the precise meaning). This will be called the "direct
method". Of course, this may be interesting whenever the expression
one obtains can be handled, which is the actual situation when the
random field has a law which is stationary and isotropic. Our method
relies on the application of some known results on the spectrum of
random matrices.
\medskip

{\bf Second}, one can use Theorem \ref{exact} to study the
asymptotics of $\P\{M>u\}$  as $ u \to +\infty$. More precisely, one
wants to write, whenever it is possible
\begin{equation} \label{genasint1}
\P\{M>u\}=~A(u)~\exp \big(-\frac{1}{2} \frac{u^2}{\sigma ^2}
\big)~+~B(u)
\end{equation}
where $A(u)$ is a known function having polynomially bounded growth
as $u\rightarrow +\infty$, $\sigma^2=\sup_{t \in S}\Var(X(t))$ and
$B(u)$ is an error bounded by a centered Gaussian density with
variance $\sigma^2_1$, $\sigma^2_1 < \sigma^2.$ We will call the
first (respectively the second) term in the right-hand side of
(\ref{genasint1}) the "first (resp second) order approximation of
$\P\{M>u \}.$"\\

First order approximation has been considered in \cite{AT}
\cite{TTA} by means of the expectation of the EPC of the excursion
set $ E_u := \{ t \in S: X(t) >u\}$. This works for large values of
$u$. The same authors have considered the second order
approximation, that is, how fast does the difference between
$\P\{M>u\}$ and the expected EPC tend to zero when $u \to +\infty$.
\medskip

We will address the same question both for the direct method and the
EPC approximation method. Our results on the second order
approximation only speak about the size of the variance of the
Gaussian bound. More precise results are only known to the authors
in the special case where $S$ is a compact interval of the real
line, the Gaussian process $ \mathcal{X}$ is stationary and
satisfies a certain number of additional requirements (see Piterbarg
\cite{piter96} and Aza\"{\i}s et al. \cite{ABW}).\\

Theorem \ref{varpasconst} is our first result in this direction.  It
gives a rough bound for the error $B(u)$ as $u\rightarrow +\infty$,
in the case the maximum variance is attained at some strict subset
of the face in $S$ having the largest dimension. We are not aware of
the existence of other known results under similar
conditions.\medskip

 In Theorem
\ref{varconst} we consider processes with constant variance. This is
close to Theorem 4.3 in \cite{TTA}. Notice that Theorem
\ref{varconst} has some interest only in case $\sup_{t\in S}\kappa_t
<\infty$, that is, when one can assure that $\sigma^2_1 < \sigma^2$
in (\ref{genasint1}). This is the reason for the introduction of the
additional hypothesis $\kappa (S)<\infty$ on the geometry of $S$,
(see below (\ref{capas}) for the definition of $\kappa (S)$), which
is verified in some relevant situations (see the discussion  before
the statement of Theorem \ref{varconst}).

\medskip

In Theorem \ref{eqisot}, $S$ is convex and the process stationary
and isotropic. We compute the exact asymptotic rate for the second
order approximation as $u\rightarrow +\infty$ corresponding to the
direct
method.\\

In all cases, the second order approximation for the direct method
provides an upper bound for the one arising from the EPC
method.\\

Our proofs use almost no differential geometry, except for
some elementary notions in Euclidean space. Let us remark also that
we have separated the conditions on the law of the process from
the conditions on the geometry of the parameter set.\\

\medskip

{\bf Third}, Theorem \ref{exact} and related results in this paper,
in fact refer  to the density  $p_M$ of the maximum. On integration,
they imply immediately a certain number of properties of the
probability distribution $F_M$, such as the behaviour of the tail as
$u \to +\infty$.
\medskip

Theorem \ref{exact} implies that  $F_M$ has a density and we have an
implicit expression for it. The proof of this fact here appears to
be simpler than previous ones (see  Aza\"{\i}s and Wschebor
\cite{AW}) even in the case the process has 1-dimensional parameter
(Aza\"{\i}s and Wschebor \cite{AW0}). Let us remark that Theorem
\ref{exact} holds true for non-Gaussian processes under appropriate
conditions allowing to apply Rice formula.\medskip

Our method can be exploited to study higher order differentiability of
$F_M$ (as it has been done in \cite{AW0} for one-parameter processes) but we will not pursue
this subject here. \medskip

This paper is organized as follows: \medskip

Section 2 includes an extension of Rice Formula  which gives  an
integral expression for the expectation  of the weighted number of
roots of a random system of $d$ equations with $d$ real unknowns. A
complete proof of this formula in a form which is adapted to our
needs in this paper, can be found in \cite{AW06}. There is an
extensive literature on Rice formula in various contexts (see for
example Belayiev \cite{belyaiev} , Cram\'er-Leadbetter
\cite{cramer},
 Marcus \cite{marcus0}, Adler \cite{A81}, Wschebor
 \cite{wschebor}.\\

In Section 3, we obtain the exact expression for the distribution of
the maximum as a consequence of the Rice-like formula of the
previous section. This immediately implies the existence  of the
density and gives the implicit formula for it. The proof avoids
unnecessary technicalities  that we have used in previous work, even
in cases that are much simpler than the ones considered here.\\

In Section 4, we  compute  (Theorem \ref{tpbar}) the first order
approximation in the direct method for stationary isotropic
processes defined on a polyhedron, from which a new upper bound for
$\P\{M>u\}$ for all real $u$ follows.
\medskip

In Section 5, we consider second order approximation, both for the
direct method and the EPC approximation method. This is the content
of Theorems \ref{varpasconst}, \ref{varconst} and \ref{eqisot}.
\medskip

Section 6 contains some examples.\medskip

\subsection*{Assumptions and notations}

$\mathcal{X}= \{X(t): t\in S\} $ denotes a real-valued   Gaussian
field  defined  on the parameter set $S$. We assume that $S$
satisfies the hypothesis  A1
\begin{itemize}
 \item[A1] :
  \item  $S$ is a compact subset of $\R^d$
  \item  $S$ is the disjoint  union of $S_d , S_{d-1} ... , S_0$,
where $S_j$ is an orientable $C^3$ manifold of dimension $j$ without
boundary. The $S_j$'s will be called faces. Let $S_{d_0}$,  $d_0
\leq d$  be the  non empty  face having largest dimension.
 \item We will assume that each $S_j$ has an atlas such that the second derivatives of the
 inverse functions of all charts (viewed as diffeomorphisms from an open set in $\R^j$ to $S_j$)
 are bounded by a fixed constant.
For $t\in S_j$, we denote $L_t$ the maximum curvature of $S_j$ at
the point $t$. It follows that $L_t$ is bounded for $t \in S$.
\end{itemize}

Notice that the decomposition $S = S_d \cup...\cup S_0$ is not
unique.

\medskip

Concerning the random field  we make the following assumptions
A2-A5
\begin{itemize}
\item[A2] : $\mathcal{ X }$ is in fact defined on an open  set containing $S$ and has $ { \cal C}^2$ paths
 \item [A3] : for every $t\in S$
   the distribution of $ \big(X(t),X'(t)\big) $ does not degenerate; for every $s,t \in S$,
$s \neq t$, the distribution of $ \big(X(s),X(t)\big) $ does not degenerate.
  \item [A4] : Almost surely the maximum of $X(t)$ on $S$ is attained at a single point.
  \end{itemize}
  For $t \in S_j$, $X'_j(t)$ $X'_{j,N}(t)$ denote respectively  the derivative along
  $ S_j$ and the normal derivative. Both quantities are viewed as vectors
   in $\R^d$,
and the density of their distribution  will be expressed
respectively with respect to an orthonormal
   basis of the tangent space $T_{t,j}$ of $ S_j$ at the point
   $t$, or its orthogonal complement $N_{t,j}$.
$X''_j(t)$ will denote  the second derivative of $X$  along $S_j$,
at the point
   $t \in S_j $ and will be viewed as a matrix expressed in an
    orthogonal basis of $T_{t,j}$. Similar notations will be used
    for any function defined on $S_j$.
\begin{itemize}
  \item [A5]  : Almost surely, for every $j=0,1,\ldots,d$  there is no point $t$ in $S_j$
  such that $ X'_j(t) =0$, $\det(X''_j(t)) =0$
\end{itemize}

Other notations and conventions will be as follows :

\begin{itemize}
  \item  $ \sigma_j$ is the geometric measure on $S_j$.
\item $ m(t):= \E(X(t))$,  $r(s,t) = \Cov(X(s),X(t))$ denote respectively the expectation and covariance of
the process
$\mathcal{X}$ ; $r_{0,1}(s,t)$, $r_{0,2}(s,t)$ are the first and
the second derivatives of $r$ with respect to  $t$. Analogous
notations will be used for other derivatives  without further
reference.
\item  If $\eta$ is a random variable  taking values in some
Euclidean space, $ p_\eta(x)$ will denote the density  of its
probability distribution with respect to the Lebesgue  measure,
whenever it exists.
\item $\varphi(x) = (2\pi)^{-1/2} \exp( -x^2/2)$ is the standard Gaussian
density ; $ \Phi(x) := \int_{-\infty}^x \varphi(y) dy$.
 \item Assume that the random vectors $ \xi,\eta$ have a joint
 Gaussian distribution, where $\eta$ has values  in some finite
 dimensional Euclidean space. When it is well defined,
 $$
 \E(f(\xi)/ \eta=x)
 $$
 is the version of the conditional expectation obtained using
 Gaussian regression.
 \item $E_u := \{t \in S : X(t)>u \}$  is the excursion set
 above $u$ of the function $ X(.)$ and $A_u := \{M \leq u\}$ is the event that the maximum
is not larger than $u$.

 \item $ \langle, \rangle, \| \| ,$ denote
 respectively  inner product and norm in a finite-dimensional
 real Euclidean space;  $ \lambda_d$ is the Lebesgue measure  on $\R^d$; $ \mathcal{S}^{d-1}$ is the unit sphere
 ; $A^c$  is the complement of the set $A$.  If $M$ is a real square matrix,
$M\succ 0$ denotes that it is positive definite.

 \item If $g: D \to C$ is a function  and $u \in C$, we denote
 $$
 N_u^g(D) := \sharp\{ t \in D: g(t) =u \}
 $$
which may be finite or infinite.

\end{itemize}

\subsection*{Some remarks on the hypotheses}

One can give simple sufficient additional conditions on the
process $ \mathcal{X}$ so that  A4 and A5 hold true.

If we assume that for each pair $j,k= 0,\ldots,d$ and each pair of
 distinct points $s,t$, $s\in S_j, t\in S_k$, the distribution of
 the triplet
 $$
 \big( X(t) -X(s) , X'_j(s), X'_k(t)) \big)
 $$
 does not degenerate in $ \R \times \R^j \times \R^k$, then  A4  holds true.
 \medskip

 This is well-known and follows easily from the next lemma (called
 Bulinskaya 's lemma) that we state without proof, for
 completeness.

\begin{lem} \label{buli} Let $Z (t)$  be  a stochastic  process defined on some
neighborhood  of a set $T$ embedded in some Euclidean space.
Assume that the  Hausdorff dimension  of $T$ is smaller or equal
than the integer  $m$ and that the values of $Z$ lie in
  $\R^{m+k}$ for some positive integer $k$ .
 Suppose, in addition,  that  $Z$  has $ \c1 $ paths  and that the density $ p_{Z(t)} (v) $ is bounded
  for $t\in T$
  and  $v$ in some neighborhood of
$u \in \R^{m+k}$. Then,  a. s.  there is no point $t \in T $  such that $
Z(t) =u$.
\end{lem}

With respect to A5, one has the following sufficient conditions:
Assume  A1, A2,  A3  and as additional hypotheses one of the
following two:
\begin{itemize}
  \item$t\Fle X(t)$  is of class $\mathcal{C}^{3}$ \\
  \item$$
\sup_{t\in S,x^{\prime }\in
V(0)}\P\big(|\det\big(X^{\prime \prime }(t)\big)|<\delta /X^{\prime
}(t)=x^{\prime }\big)\rightarrow 0,~~\mbox{ as } \delta \to 0,
$$
where $ V(0)$ is some neighborhood of zero.
\end{itemize}
Then A5 holds true. This follows  from  Proposition 2.1 of
\cite{AW} and \cite{CW}.

\section{Rice formula for the number of  weighted roots of
random fields} \label{srice}
 In this section we review Rice formula for the expectation of the
 number of roots of a random system of equations. For proofs,
 see for example \cite{AW}, or  \cite{AW06}, where a simpler one is given.\medskip

\begin{theo}[Rice formula] \label{sanspoids}

Let $Z : U \to \R^d$ be a random field, $U$ an open subset of
$\R^d$  and $u \in \R^d$ a fixed point  in the codomain. Assume that: \\
(i) $Z$ is Gaussian,~~ \\
(ii) almost surely  the function  $t\rightsquigarrow Z(t)$ is of
class $\mathcal{C}^{1}$, \\
(iii)  for each $t\in U$, $Z(t)$ has a non degenerate distribution
(i.e. $\Var \big(Z(t)\big) \succ 0)$, \\
(iv) $\P\{\exists t \in U,Z(t)=u,\det\big(Z^{\prime }(t)\big)=0\}=0$

Then, for every Borel set B  contained in $U$,  one has
\begin{equation} \label{rice}
\E \left( N_{u}^{Z}(B)\right) =\int_{B}\E\left( |\det (Z^{\prime
}(t))|/Z(t)=u\right) ~p_{Z(t)}(u)dt.
\end{equation}
If B is compact, then both sides in (\ref{rice}) are finite.
\end{theo}

\begin{theo} \label{poids}
 Let $ Z$ be a random field that verifies
the hypotheses of Theorem \ref{sanspoids}. Assume that for each
$t\in U $  one has another random field $Y^t : W \to \R^{d'}$, where
$W$ is some topological space, verifying the following conditions:
\begin{itemize}
  \item [a)] $Y^t(w)$ is a measurable function of $( \omega,t,w)$
  and almost surely, $(t,w)\rightsquigarrow Y^t(w)$ is
  continuous.

  \item[b)] For each $t\in U$  the random process $(s,w) \Fle \big(
  Z(s), Y^t(w) \big)$  defined on $U \times W$  is Gaussian.
\end{itemize}
Moreover, assume that $g: U \times \mathcal{C}(W, \R^{d'}) \to \R $
is a bounded function, which is continuous when one puts on
$\mathcal{C}(W, \R^{d'})$ the topology of uniform convergence on
compact sets. Then, for each compact subset $I$ of $U$, one has
\begin{equation}\label{forpoids}
   \E\big( \sum_{t\in I, Z(t) =u} g(t,Y^t) \big)
   =  \int_I \E\big( |\det(Z'(t)| g(t,Y^t)/Z(t) =u) . p_{Z(t)} (u)
   dt.
\end{equation}
\end{theo}
{\bf Remarks:}

1. We have already mentioned in the previous section sufficient
conditions implying hypothesis (iv) in Theorem \ref{sanspoids}.

2. With the hypotheses of Theorem \ref{sanspoids} it follows
easily that if $J$ is a subset of $U$, $\lambda_d(J)=0$, then $\P\{ N_u^Z(J) =0\}
=1$ for each $u \in \R^d$.\medskip

\section{The implicit formula for the density of the maximum}

\begin{theo}\label{exact}
Under assumptions A1 to A5, the distribution of $M$ has the density
\begin{align}
p_M(x) &=  \sum _{t \in S_0} \E\big( \UN _{A_x} / X(t)=x \big)
p_{X(t)} (x) \notag
\\
\label{densm}
  &+
 \sum _{j=1}^d \int_{S_j}
\E\big(|\det(X''_j(t))| \UN _{A_x} /X(t)=x, X'_j(t) =0\big) p_{X(t),
X'_j(t)} (x,0) \sigma_j(dt)  ,
\end{align}
\end{theo}

{\bf Remark:}  One can replace $|\det(X''_j(t))| $ in
the conditional expectation  by $ (-1)^j \det(X''_j(t))$, since
under the conditioning and whenever $M\leq x $ holds true,
$X''_j(t)$ is negative semi-definite.

\subsubsection*{ Proof of Theorem \ref{exact} } Let
$N_j(u),j=0,\ldots,d$ be the number of global maxima  of $X(.)$ on
$S$ that belong to $S_j$ and are larger than $u$. From the
hypotheses it follows that a.s.  $ \sum_{j=0,\ldots,d} N_j(u)$ is
equal  to 0 or 1, so that
\begin{equation}\label{sumn}
  \P\{M>u\}  =\sum_{j=0,\ldots,d} \P\{N_j(u) =1 \} =\sum_{j=0,\ldots,d}
  \E(N_j(u)).
\end{equation}
The proof will be finished  as soon as we show that each term in
(\ref{sumn})  is the integral over $ (u,+\infty)$  of the
corresponding term in (\ref{densm}). \\
This is self-evident  for $j=0$. Let us consider the term $j=d$.
We apply the weighted Rice formula of Section \ref{srice} as
follows :
\begin{itemize}
  \item $Z$ is the random field $X'$  defined on  $S_d$.

  \item  For each $t\in S_d $, put $W=S $ and $Y^t : S \to \R^2 $
  defined as:
  $$
  Y^t(w):= \big( X(w) -X(t) ,X(t) \big).
  $$
  Notice that the second coordinate in the definition of $Y^t$ does
  not depend on $w$.
  \item  In the place of the function $g$, we take for each $n
  = 1,2,\ldots$ the function $g_n$  defined as follows:
  $$
  g_n (t,f_1,f_2)=~g_n (f_1,f_2) = \big( 1- \mathcal{ F}_n (\sup_{ w \in S}
  f_1(w)) \big) . \big( 1- \mathcal{ F}_n ( u-f_2(\overline{w}))\big),
  $$
  where $\overline{w}$ is any point in $W$ and for $n $ a positive
 integer and $x \geq 0$, we define :
\begin{equation}\label{fm}
  \mathcal{F}_n (x)  := \mathcal{F} (nx) ~~;~~\mbox{ with }
  \mathcal{F}(x)= 0 \mbox{ if } 0 \leq
x\leq 1/2~~,~~\mathcal{F}(x) =1 \mbox{ if }  x \geq 1~,
\end{equation}
and $\mathcal{F}$ monotone  non-decreasing
and continuous.
\end{itemize}
It is easy  to check that all the requirements in Theorem
\ref{poids} are satisfied, so that, for the value 0 instead of
$u$ in formula (\ref{forpoids}) we get:
\begin{equation}\label{forpoidsn}
   \E\big( \sum_{t\in S_d, X'(t) =0} g_n(Y^t) \big)
   =  \int_{S_d}\E\big( |\det(X''(t)| g_n(Y^t)/X'(t) =0) . p_ {X'(t)}  (0)
  \lambda_d(dt).
\end{equation}
Notice that the formula holds true for each compact subset of $S_d$
in the place of $S_d$, hence for $S_d$ itself by monotone convergence.\\
Let now $n \to \infty$  in (\ref{forpoidsn}). Clearly  $g_n(Y^t)
\downarrow  \UN_{X(s)-X(t) \leq 0, \forall s \in S}. \UN _{X(t)
\geq u} $.  The passage to the limit does not present any
difficulty since $ 0 \leq g_n(Y^t) \leq 1 $ and the sum   in the
left-hand side is bounded by the random variable
$N_0^{X'}(\overline{S_d})$, which is in $L^1$ because of Rice
Formula. We get
$$
 \E(N_d  (u) ) = \int_{S_d}\E\big( |\det(X''(t)|  \UN_{X(s)-X(t) \leq 0, \forall s \in S} \UN _{X(t)
\geq u} /X'(t) =0) . p_ {X'(t)}  (0)
  \lambda_d(dt)
  $$
  Conditioning on the value of $X(t)$, we obtain the desired
  formula for $j=d$. \medskip

The proof for $1 \leq j \leq d-1$  is essentially the same, but
one must take care  of the parameterization of the manifold $S_j$.
One can first establish locally the formula on a chart of $S_j$,
using local coordinates.

It can be proved  as  in \cite{AW},  Proposition 2.2 (the only
modification is due to the term $\UN _{A_x } $)
 that the quantity written in some chart as
$$
\E\big(\det(Y''(s))
\UN _{ A_x} /Y(s) =x, Y'(s) =0\big)
p_{Y(s),Y'_j(s)} (x,0) ds,
$$
where the process $Y(s)$ is the process $X$ written in some chart $ $  of $S_j$ ,\\
($Y(s) = X( \phi^{-1}(s))$),
 defines a $j$-form. By a $j$-form we mean  a mesure on $S_j$ that does not depend on the parameterization and
   which has a density with respect to the Lebesgue measure $ds$ in every chart. It can be  proved also
   that the
   integral of this $j$-form on $ S_j$ gives the expectation of $N_j(u)$.

To get formula (\ref{rice}) it suffices to consider locally around
a precise point $t \in S_j$
 the chart $ \phi$  given by the projection on the tangent space at $t$. In this  case we obtain that at $t$
\begin{itemize}
  \item  $ds$ is in fact $\sigma_j (dt)$
  \item  $Y'(s)$ is isometric to  $X'_j(t)$
\end{itemize}
where $s = \phi(t)$.\hfill $\Box $\\
\medskip

The first consequence of Theorem \ref{exact}  is the next
corollary. For the statement, we need  to introduce some further
notations.

For $t$ in $ S_j$, $ j\leq d_0$ we define $\C $ as the closed convex
cone generated by the set of directions:
 $$
  \{\lambda  \in \R^d :~\|\lambda\|=1~; \exists ~s_n \in
  S, (n=1,2,\ldots) \mbox{ such that }
  s_n \to t, \frac{t-s_n }{\|t-s_n
 \|}\to \lambda \mbox{ as } n \to + \infty \},
 $$
whenever this set is non-empty and $ \C  =\{0\}$ if it is empty. We
will denote by $\ch $ the dual cone of $ \C $, that is:
$$
\ch := \{z~\in~\R^d : \langle z , \lambda \rangle \geq 0 \mbox{ for
all } \lambda \in \C  \}.
$$
Notice that these definitions easily imply that $T_{t,j} \subset
C_{t,j}$ and $\ch  \subset N_{t,j}$. Remark also that for $j=d_0$,
$\ch =N_{t,j}$.

We will say that the  function $X(.)$ has an "extended outward"
derivative at the point $t$ in $S_j$, $ j \leq d_0$   if
$X'_{j,N}(t) \in \ch $.

\begin{cor} \label{cpbar} Under assumptions A1 to A5, one has :
\begin{itemize}
  \item [(a)] $ p_M(x) \leq \overline{ p}(x)$ where
\begin{multline}\label{densmaj}
 \overline{ p}(x):=\sum _{t \in S_0}  \E\big ( \UN_{X'(t)  \in
 \widehat{C}_{t,0}} / X(t) =x\big) p_{X(t)} (x) +
 \\
\sum _{j=1}^{d_0} \int_{S_j}
 \E\big(|\det(X''_j(t))| \UN_{X'_{j,N}(t)\in \ch}
/ X(t) =x, X'_j(t)=0\big) p_{X(t), X'_j(t)} (x,0) \sigma_j(dt) .
\end{multline}

\item [(b)]$ \P\{ M > u\} \leq \displaystyle \int _u ^{+\infty} \overline{ p}(x) dx.$
\end{itemize}
 \end{cor}

 %
 %
\subsubsection*{Proof} (a) follows from Theorem \ref{exact} and the
observation that if $ t\in S_j,$ one has \\
$ \{M \leq X(t) \} \subset \{ X'_{j,N}(t) \in \ch \}$. (b) is an obvious
consequence of (a).\hfill $\Box $\\

The actual interest of this Corollary depends on the feasibility of
computing $\overline{p}(x)$. It turns out that it can be done in
some relevant cases, as we will see in the remaining of this
section. Our result can be compared with the approximation of $\P\{
M > u\}$ by means of $ \int _u^{+\infty} p^E(x) dx  $ given by
\cite{AT}, \cite{TTA} where
\begin{multline}\label{denseu}
p^E(x):=\sum _{t \in S_0}  \E\big ( \UN_{X'(t)  \in
 \widehat{C}_{t,0}} / X(t) =x\big) p_{X(t)} (x)\\
 +
\sum _{j=1}^{d_0} (-1)^j  \int_{S_j}
 \E\big(\det(X''_j(t)) \UN_{X'_{j,N}(t)\in \ch}
/ X(t) =x, X'_j(t) =0\big) p_{X(t), X'_j(t)} (x,0) \sigma_j(dt) .
\end{multline}
Under certain conditions ,
 $\int_u ^{+\infty}p^E(x)  dx $ is the expected value of the
 EPC  of the excursion set $E_u$ (see \cite{AT}).
 The advantage  of $p^E(x)$ over $ \overline{p}(x)$  is that one
 can have nice expressions for it in quite general situations.
 Conversely $ \overline{p}(x)$ has the obvious advantage that it
 is an upper-bound of the true density $ p_M(x)$ and hence
 provides upon integrating once, an upper-bound for the tail probability,
 {\bf for every $u$ value}. It is not
 known whether a similar inequality holds true for  $p^E(x)$.
 \medskip \\
 On the other hand, under additional conditions, both provide good first order approximations
 for $ p_M(x)$  as $ x \to \infty$ as we will see in the next
 section. In the special case in which the process $\mathcal{ X} $ is centered
 and has a law that is invariant under isometries and translations, we describe
  below a procedure to compute $\overline{ p}(x)$.

 \section{ Computing $\overline{ p}(x)$ for stationary
 isotropic Gaussian fields } \label{comppbar}

 For one-parameter centered Gaussian process having constant
 variance and satisfying certain regularity conditions, a general
 bound for $ p_M(x)$ has been computed in \cite{AW}, pp.75-77. In the two
 parameter case, Mercadier \cite{mercadier} has shown a
 bound for $ \P\{M>u\} $, obtained by means of a method
 especially suited to dimension 2. When the parameter is one or two-dimensional,
 these bounds are sharper
 than the ones below which, on the other hand, apply to any
 dimension but to a more restricted context.
We will assume now that the process $\mathcal{ X}$ is centered
Gaussian,  with a  covariance  function that can be written as
\begin{equation}\label{isotro}
  \E\big( X(s).X(t) \big)  = \rho\big( \| s-t\|^2\big),
\end{equation}
where $ \rho: \R^+ \to \R $ is of class $ \mathcal{C}^4$ . Without
loss of generality, we assume that $ \rho(0) =1$.  Assumption
(\ref{isotro}) is equivalent  to saying that the law of $\mathcal{
X}$ is invariant under isometries (i.e. linear transformations that
preserve the scalar product) and translations of the underlying parameter space $\R^d$.

We will also assume that the set $S$ is a polyhedron. More precisely
we assume that each $S_j  (j= 1,\ldots,d)$ is a union  of subsets of
affine  manifolds of dimension $j$ in $\R^d$.
\medskip

The next lemma  contains some auxiliary computations which are
elementary and left to the reader. We use the abridged notation :
$\rho':=\rho'(0)$, $\rho'':=\rho''(0)$
%
\begin{lem}\label{formules}
 Under the conditions above, for each $t \in U$, $ i,i',k,k',j =1,\ldots ,d$:
\begin{enumerate}
  \item $ \E \big( \frac{\partial X}{\partial t_i}(t). X(t)\big)  =0$,
  \item $ \E \big( \frac{\partial X}{\partial t_i}(t). \frac{\partial X}{\partial t_k}(t)\big)  =
  - 2 \rho' \delta_{ik}
  \mbox { and }\rho' <0$,
  \item $ \E \big( \frac{\partial ^2 X}{\partial t_i \partial t_k}(t). X(t)\big)  =
  2 \rho' \delta_{ik},
\E \big( \frac{\partial ^2 X}{\partial t_i \partial t_k}(t). \frac{\partial X}{\partial t_j}(t)\big)  = 0 $
 \item $ \E \big( \frac{\partial ^2 X}{\partial t_i \partial t_k}(t).
 \frac{\partial ^2 X}{\partial t_{i'} \partial t_{k'}}(t)
 \big)  = 24 \rho''
 \big[  \delta_{ii'}.\delta_{kk'}+ \delta_{i'k} .\delta_{ik'} +\delta_{ik}
 \delta_{i'k'}$\big],
\item $ \rho''- \rho'^2 \geq 0$
  \item  If $ t \in S_j$,  the conditional distribution  of $ X_j^{\prime
  \prime}(t)$ given $ X(t) = x , X'_j(t) =0$ is the same as the
  unconditional distribution of the random matrix
  $$
  Z + 2 \rho' x I_j~,
$$
where  $ Z= (Z _{ik}: i,k= 1,\ldots,j)$ is a symmetric $j \times j $
matrix  with centered Gaussian entries, independent  of the pair $
\big( X(t), X'(t)\big)$ such that, for $ i \leq k$, $i'\leq k'$
one has :
$$
\E(Z_{ik} Z_{i'k'}) = 4 \big [ 2 \rho'' \delta_{ ii'} + (
\rho''-\rho'^2) \big] \delta_{ik}
 \delta_{i'k'} +4 \rho'' \delta_{ii'}.\delta_{kk'}(1-\delta_{ik}) ~.
 $$
\end{enumerate}
\end{lem} \medskip

Let us introduce some additional notations:
\begin{itemize}
  \item $H_n(x), n= 0,1,\ldots $  are the standard Hermite
  polynomials, i.e.
  $$
  H_n(x) := e^{ x^2} \big( - \frac{\partial}{\partial x} \big)^n e^{ -x^2}
  $$
  For the properties of the Hermite polynomials we refer to Mehta
  \cite{mehta}.
\item $\overline{H}_n(x), n= 0,1,\ldots $  are the modified  Hermite
  polynomials, defined as:
  $$
 \overline{ H}_n(x) := e^{ x^2/2} \big( - \frac{\partial}{\partial x} \big)^n e^{ -x^2/2}
  $$
\end{itemize}
We will use the following result:
 \begin{lem} \label{lemjn}
 Let
\begin{equation}\label{jotan}
J_n(x) := \int_{-\infty }^{+\infty} e^{ -y^2/2} H_n(\nu)
dy,~~n=0,1,2,\ldots
\end{equation}
where $\nu$ stands for the linear form $ \nu = a y + bx$ where
$a,b$ are some real parameters that satisfy $a^2+b^2= 1/2$. Then
$$
J_n(x) := (2b)^n \sqrt{2 \pi}~\overline{ H}_n(x).
$$
\end{lem}
\subsubsection*{Proof :} It is clear that $J_n$ is a polynomial
having degree $n$. Differentiating in (\ref{jotan}) under the
integral sign, we get:
\begin{equation}\label{jota2}
J'_n(x) = b \int_{-\infty }^{+\infty} e^{ -y^2/2} H'_n(\nu) dy
=2nb \int_{-\infty }^{+\infty} e^{ -y^2/2} H_{n-1}(\nu) dy = 2 n ~b~
J_{n-1} (x)
\end{equation}
Also:
$$
J_n(0) = \int_{-\infty }^{+\infty} e^{ -y^2/2} H_n(ay) dy,
$$
so that $J_n(0) = 0$ if $n$ is odd.\\
If $n$ is even, $n \geq 2$, using the standard recurrence
relations for Hermite polynomials, we have:
\begin{align} \label{jota3}
J_n(0) &= \int_{-\infty }^{+\infty} e^{ -y^2/2} \big[ 2ay
H_{n-1}(ay)  -2(n-1) H_{n-2}(ay)\big] dy \notag
\\
&= 2 a^2 \int_{-\infty }^{+\infty} e^{ -y^2/2} H'_{n-1}(ay) dy -
2(n-1) J_{n-2} (0) \notag
\\
&= -4 b^2(n-1) J_{n-2}(0).
\end{align}
Equality (\ref{jota3}) plus $J_0(x) = \sqrt{2 \pi}$ for all $x \in\R$,
imply that:
\begin{equation}\label{jota4}
  J_{2p} (0) = (-1)^p (2b)^{2p} (2p-1)!! \sqrt{2 \pi} =(-2b^2)^p
\frac{  (2p)!}{p!}~\sqrt{2 \pi}.
\end{equation}
Now we can go back to (\ref{jota2}) and integrate successively for
$n=1,2,\ldots$ on the interval $[0,x]$ using the initial value
given by (\ref{jota4}) when $n=2p$ and $J_n(0) =0$ when $n$ is
odd, obtaining :
$$
J_n(x) = (2b)^n \sqrt{2 \pi} Q_n(x) ,
$$
where the sequence of polynomials $ Q_n, n=0,1,2,\ldots$ verifies
the conditions:
\begin{align}
Q_0(x) &= 1 \\
Q'_n(x) &= nQ_n(x) \\
 Q_n(0) &= 0 ~\mbox{ if } n \mbox{ is odd }  \\
Q_n(0) &= (-1)^{n/2}(n-1)!! ~\mbox{ if } n \mbox{ is even. }
\end{align}
It is now easy to show that in fact $Q_n(x) =
\overline{H}_n(x)~,~n=0,1,2,\ldots$ using for example  that:
$$
\overline{H}_n(x) = 2^{n/2} H_n \big( \frac{x}{\sqrt{2}}\big).
$$
\hfill $\Box $\\

The integrals
$$
I_n(v) = \int_{v}^{+\infty} e^{ -t^2/2} H_n(t) dt,
$$
will appear in our computations. They are computed in the next Lemma, which can be
proved easily, using the standard properties of Hermite polynomials.

\begin{lem} \label{lemin}
(a)\begin{align}\label{in}
I_n(v) & =~2e^{ -v^2/2}\sum _{k=0}^{[\frac{n-1}{2}]}2^k~\frac{(n-1)!!}{(n-1-2k)!!}H_{n-1-2k}(v)
\\
       & +\UN_{ \{ n\mbox{\rm { even}}\} }~2^{\frac{n}{2}}~(n-1)!!\sqrt{2\pi}~\overline{\Phi}(x)
       \end{align}

(b)
\begin{equation}\label{in2}
I_n(-\infty) =
       \UN_{ \{n\mbox{\rm { even}}\}}2^{\frac{n}{2}}~(n-1)!!~\sqrt{2\pi}
\end{equation}
\end{lem}

\begin{theo} \label{tpbar}Assume that the process $\mathcal{X}$ is centered Gaussian,
satisfies  conditions A1-A5 with a covariance having the form
(\ref{isotro}) and verifying the regularity conditions of the
beginning of this section. Moreover, let $S$ be a polyhedron.  Then,
$ \overline{ p}(x)$ can be expressed by means of the following
formula:
\begin{equation}\label{pbar}
  \overline{ p}(x) = \varphi(x) \left\{ \sum_{t \in S_0}
  \widehat{\sigma}_{0}(t)
   + \sum _{j=1}^{d_0} \big[ \big(\frac{|\rho'|}{\pi}\big)^{ j/2}
   \overline{H}_j(x) + R_j(x)\big] g_j\right\}~,
\end{equation}
where
\begin{itemize}
\item  $g_j $ is a geometric parameter  of the face $S_j$ defined
by
\begin {equation}\label{geom}
g_j  = \int _{S_j}  \widehat{\sigma}_{j}(t)  \sigma_{j} (dt),
\end{equation}
where $\widehat{\sigma}_{j}(t)$ is the normalized solid angle  of the cone $\ch$ in
$ N_{t,j}$, that is:
\begin{align}
 \widehat{\sigma}_{j}(t) &= \frac{\sigma _{d-j-1} ( \ch \cap \mathcal{S}^{d-j-1})}{\sigma _{d-j-1} ( \mathcal{ S}^{d-j-1})}
 \mbox{ for } j= 0,\ldots,d-1,
 \\
\widehat{\sigma}_{d}(t) &= 1.
\end{align}
Notice that for convex or other usual polyhedra
$\widehat{\sigma}_{j}(t)$ is constant for $ t \in S_j$, so that
$g_j$ is equal to this constant multiplied by the $j$-dimensional
geometric measure of $S_j$.

\item For $j= 1,\ldots d $,
\begin{equation}\label{errj}
R_j (x) = \big(\frac{2\rho''}{\pi|\rho'|}\big)^{\frac{j}{2}} \frac{\Gamma((j+1)/2}{\pi}
\int_{-\infty}^{+\infty}~T_j(v)\exp\big(- \frac{y^2}{2}\big)  \ dy
\end{equation}
\\
 where
\begin{equation}\label{nu}
 v := - (2)^{ -1/2} \big( (1- \gamma^2) ^{1/2} y - \gamma x \big)
  ~\mbox{ with }~ \gamma := |\rho'| (\rho'')^{ -1/2}
  \end{equation}
and
\begin{equation}\label{tejota}
 T_j(v):= \big[\sum _{k=0}^{j-1}\frac{H_k^2(v)}{2^kk!}\big]e^{ -v^2/2}~-
~\frac{H_j(v)}{2^j(j-1)!}I_{j-1}(v).
  \end{equation}
\\
where $I_n$ is given in the previous Lemma.
\end{itemize}
\end{theo}

For the proof of the theorem, we need some ingredients from random
matrices theory. Following Mehta \cite{mehta}, denote by $q_n(\nu )$
the density  of eigenvalues of $n \times n $ GOE matrices at the
point $\nu$, that is, $q_n(\nu) d\nu$  is the probability of $G_n$
having an eigenvalue in the interval $ (\nu, \nu+d\nu) $.
The random $n \times n$ real random matrix  $G_n$  is said to have the  GOE
distribution, if it is symmetric, with centered
Gaussian  entries $g_{ik}, i,k=1,\ldots,n$ satisfying $
\E(g_{ii}^2 )= 1$, $\E(g_{ik}^2)= 1/2$ if $i<k$ and the
random variables: $ \{g_{ik}, ~1 \leq i \leq k \leq n\}$ are
independent.\\
   It is well known
  that:
  \begin {align}
  e^{\nu^2/2}
 q_n(\nu)  &= e^{-\nu^2/2}\sum _{k=0}^{n-1} c_k^2 H_k^2(\nu)  \notag   \\
           &+  1/2 ~(n/2)^{1/2} c_{n-1} c_n H_{n-1} (\nu)
                   \Big[ \int_{-\infty}^{+\infty}  e^{-y^2/2}
                   H_n(y) dy -2\int_{\nu}^{+\infty}e^{-y^2/2}
                   H_n(y) dy \Big] \notag   \\
 &+ \UN_{\{n \mbox{ odd } \}} \frac{H_{n-1} (\nu)}{\int_{-\infty}^{+\infty}  e^{-y^2/2}
                   H_{n-1} (y) dy  }, \label{fmehta}
\end{align}
where $c_k := (2^k k! \sqrt{\pi} )^{-1/2}, k=0,1,\ldots$, (see
Mehta \cite{mehta}, ch. 7.)

%
%
%
In the proof of the theorem we will use the following  remark  due
to Fyodorov \cite{fyodorov} that we state as a Lemma
\begin{lem}\label{lfyod}
 Let $G_n$ be a GOE $n \times n$ matrix.
Then, for $ \nu \in \R$ one has:
\begin{equation}\label{fyo}
\E\big( |\det(G_n-\nu I_n)|\big) = 2^{3/2} \Gamma\big((n+3)/2\big)
\exp( \nu^2/2)\frac{ q_{n+1}(\nu)}{n+1} ,
\end{equation}
\end{lem}

\subsubsection*{Proof:}
Denote by $ \nu_1,\ldots,\nu_n$ the eigenvalues of $G_n$. It is
well-known (Mehta \cite{mehta}, Kendall et al. \cite{kendall})
that the joint density $f_n$ of the $n$-tuple of random variables
$( \nu_1,\ldots,\nu_n)$ is given by the formula
$$
f_n( \nu_1,\ldots,\nu_n) = c_n \exp\Big(  - \frac{\sum_{i=1} ^n
\nu_i^2}{2}\Big) \prod_{ 1 \leq i < k \leq n }| \nu_k - \nu_i |~~,~~
\mbox{ with } c_n := (2\pi)^{-n/2}( \Gamma(3/2)) ^n  \big(
\prod_{i=1} ^n  \Gamma(1+i/2)\big)^{-1}
$$
Then,
\begin{multline*}
 \E\big( |\det(G_n-\nu I_n)|\big)
 =\E\big(  \prod_{i=1}^n | \nu_i -\nu| \big) \\
 = \int_{\R^n} \prod _{i=1}^n | \nu_i -\nu| c_n \exp(  - \frac{\sum_{i=1} ^n
\nu_i^2}{2}) \prod_{ 1 \leq i < k \leq n} | \nu_k - \nu_i |~
d\nu_1,\ldots, d\nu_n  \\
= e^{\nu^2/2}\frac{ c_n }{c_{n+1}}\int_{\R^n}
f_{n+1}(\nu_1,\ldots,\nu_n,\nu)d\nu_1,\ldots, d\nu_n  =
e^{\nu^2/2}\frac{ c_n }{c_{n+1}} \frac{q_{n+1}(\nu)}{n+1}~.
\end{multline*}
The remainder is plain. \hfill $\Box $\

\subsubsection*{Proof of Theorem \ref{tpbar}:}

We use the definition (\ref{densmaj}) given in Corollary
\ref{cpbar} and the moment computations of Lemma \ref{formules}
which imply that:
\begin{align}
p_{X(t)} (x) &= \varphi(x)   \\
p_{X(t),X_j'(t)} (x,0) &= \varphi(x) (2\pi)^{-j/2} (-2
\rho')^{-j/2}  \label{pxx'} \\
X'(t) & \mbox{ is independent of }  X(t)   \\
X'_{j,N}(t) & \mbox{ is independent of } ( X''_j(t),X(t),X'_j(t) ).
\end{align}
Since the distribution of $X'(t)$ is centered Gaussian with
variance $ -2 \rho' I_d$, it follows that :
$$
\E(\UN_{X'(t) \in \widehat{C}_{t,0}}/ X(t) =x) = \widehat{\sigma}_{0}(t)
~~\mbox{ if } t\in S_0,
$$
and if $ t \in S_j,j \geq 1$:
\begin{multline} \label{zeta}
\E(|\det(X''_j(t))|\UN_{X'_{j,N}(t) \in \ch} /X(t) =x, X'_j(t)
=0) \\
= \widehat{\sigma}_j(t)~\E(|\det(X''_j(t))| /X(t) =x, X'_j(t)
=0)
\\
=\widehat{\sigma}_j(t)~\E(|\det(Z +2\rho' x I_j) |) .
\end{multline}
In the formula above, $\widehat{\sigma}_j(t)$ is the normalized solid
angle defined in the statement of the theorem and the random
$j\times j $ real matrix $Z$ has the distribution of Lemma
\ref{formules} . \\
A standard moment computations  shows that $Z$ has the same
distribution as the random matrix:
$$
 \sqrt{8\rho''} G_j + 2 \sqrt{ \rho''-\rho'^2} \xi I_j,
 $$
 where $ G_j$ is a $j\times j$ GOE random matrix, $ \xi$ is
 standard normal in $\R $ and independent of $G_j$.  So, for
 $j\geq 1$  one has
 $$
 \E\big( |\det (Z + 2 \rho' x I_j)|\big)  = (8 \rho'')^{j/2}
  \int_{-\infty}^{+\infty}\E\big( |\det (G_j - \nu I_j) | \big)
  \varphi(y) dy ,
  $$
  where $\nu$  is given by (\ref{nu}).\\
  For the conditional expectation  in (\ref{densmaj}) use   this last expression in (\ref{zeta})
  and (\ref{lfyod}). For the density in (\ref{densmaj}) use (\ref{pxx'}).  Then Lemma \ref{lemjn}
gives (\ref{pbar}). \hfill $\Box $\newline

\subsubsection*{ Remarks on the theorem}

\begin{itemize}
  \item The "principal term" is
\begin{equation}\label{pt }
   \varphi(x) \left\{ \sum_{t \in S_0}
  \widehat{\sigma}_{0}(t)
   + \sum _{j=1}^{d_0} \big[ \big(\frac{|\rho'|}{\pi}\big)^{ j/2}
   \overline{H}_j(x) \big] g_j\right\}~,
\end{equation}
which is the  product of a standard Gaussian density times a polynomial with degree $d_0$.
Integrating once, we get -in our special case- the  formula for the expectation of the EPC
 of the excursion set as given by \cite{AT}

  \item  The "complementary term" given by
  \begin{equation}\label{ct}
 \varphi(x) \sum_{j=1} ^{d_0}  R_j(x) g_j,
  \end{equation}

can be computed by means of a formula, as it follows from the
statement of the theorem above. These formulae will be in general
quite unpleasant due to the complicated form of $T_j(v)$. However,
for low dimensions they are simple. For example:

\begin{equation}\label{t1}
 T_1(v)=~\sqrt{2\pi}\big [\varphi(v)-v(1-\Phi (v))\big],
  \end{equation}
\begin{equation}\label{t2}
 T_2(v)=~2\sqrt{2\pi}\varphi(v),
  \end{equation}
\begin{equation}\label{t3}
 T_3(v)=~\sqrt{\frac{\pi}{2}}\big [3(2v^2+1)\varphi(v)-(2v^2-3)v (1-\Phi (v))\big].
  \end{equation}

\item Second order asymptotics for $p_M(x)$ as $x \to +\infty$ will be mainly
considered in the next section. However, we state already that the complementary term
 (\ref{ct}) is equivalent, as $x \to +\infty$, to

\begin{equation}\label{eqct}
\varphi(x)~
g_{d_0}K_{d_0}x^{2d_0-4}~e^{-\frac{1}{2}\frac{\gamma^2}{3-\gamma^2}x^2},
  \end{equation}

where the constant $K_j$, $j=1,2,...$ is given by:

\begin{equation}\label{kj}
K_j=2^{3j-2} \frac{\Gamma\big(\frac{j+1}{2}\big)}{\sqrt\pi(2\pi\gamma)^{j/2}(j-1)!}
\rho''^{j/4}\big(\frac{\gamma}{3-\gamma^2}\big)^{2j-4}.
  \end{equation}

We are not going to go through this calculation, which is elementary but requires some work.
An outline of it is the following. Replace the Hermite polynomials in the expression for $T_j(v)$ given by
(\ref{tejota}) by the well-known expansion:

\begin{equation}
H_j(v)=j!~\sum _{i=0}^{[j/2]} (-1)^i \frac{(2v)^{j-2i}}{i!(j-2i)!}
  \end{equation}

and $I_{j-1}(v)$ by means of the formula in Lemma \ref{lemin}.

Evaluating the term of highest degree in the polynomial part, this allows to prove that, as
$v \to +\infty$, $T_j(v)$ is equivalent to

\begin{equation}
\frac{2^{j-1}}{\sqrt\pi(j-1)!}~v^{2j-4}e^{-\frac{v^2}{2}}.
  \end{equation}

Using now the definition of $R_j(x)$ and changing variables in the
integral in (\ref{errj}), one gets for $R_j(x)$ the equivalent:

\begin{equation}
K_jx^{2j-4}~e^{-\frac{1}{2}\frac{\gamma^2}{3-\gamma^2}x^2}.
\end{equation}

In particular, the equivalent of (\ref{ct}) is given by the highest order non-vanishing term
in the sum.

\item Consider now the case in which $S$ is the sphere $\mathcal{S}^{d-1}$ and the process satisfies
the same
conditions as in the theorem. Even though the theorem can not be applied directly, it is possible
to deal with this example to compute $\overline{p}(x)$, only performing some minor changes. In this
case, only the term that corresponds to $j=d-1$ in (\ref{densmaj}) does not vanish,
$\widehat{C}_{t,d-1}=N_{t,d-1}$, so that $\UN_{X'_{d-1,N}(t) \in \widehat{C}_{t,d-1}}=1$
for each $t \in \mathcal{S}^{d-1}$ and one can use invariance under rotations to obtain:

\begin{equation}\label{esfera}
\overline{p}(x) = \varphi(x)~\frac{\sigma_{d-1}\big(\mathcal{S}^{d-1}\big)}{(2\pi)^{(d-1)/2}}
\E\big( |\det (Z + 2 \rho' x I_{d-1})+(2|\rho'|)^{1/2}\eta I_{d-1}|\big)
\end{equation}

where $Z$ is a $(d-1) \times (d-1)$ centered Gaussian matrix with the covariance structure of
Lemma \ref{formules} and $\eta$ is a standard Gaussian real random variable, independent
of $Z$. (\ref{esfera}) follows from the fact that the normal derivative at each point is centered
Gaussian with variance $2|\rho'|$ and independent of the tangential derivative. So, we apply the
previous computation, replacing $x$ by $x+(2|\rho'|)^{-1/2}~\eta$ and obtain the expression:
\begin{align}
\overline{p}(x) &= \varphi(x)~\frac{2\pi^{d/2}}{\Gamma(d/2)} \notag
\\
& \int_{-\infty}^{+\infty} \big[ \big(\frac{|\rho'|}{\pi}\big)^{(d-1)/2} \overline{H}_{d-1}
(x+(2|\rho'|)^{-1/2}y)+R_{d-1}(x+(2|\rho'|)^{-1/2}y) \big]\varphi(y)dy.
\end{align}
\end{itemize}


\section{ Asymptotics as $x \to +\infty$} \label{asymp}

In this section we will consider the errors in the direct and the
EPC methods for large values of the argument $x$.
  Theses errors are:
\begin{multline}\label{erdi}
\overline{p}(x)- p_M(x)  =\sum _{t \in S_0}
\E \big(\UN_{X'(t)\in \widehat{C}_{t,0}}  . \UN _{M>x}
 /  X(t) =x \big) p_{X(t)} (x)
\\
+  \sum _{j=1}^{d_0}
 \int_{S_j}
 \E\big( |\det(X''_j(t)| \UN_{X_{j,N}'(t)\in \ch} . \UN _{M>x}\big)
 /  X(t) =x,X'_j(t) =0\big) p_{X(t), X'_j(t)} (x,0) \sigma_j(dt)  .
\end{multline}
\begin{multline}\label{erepc}
p^E(x)- p_M(x)  =\sum _{t \in S_0} \E \big(\UN_{X'(t)\in
\widehat{C}_{t,0}}  . \UN _{M>x}
 /  X(t) =x \big) p_{X(t)} (x)
\\
+  \sum _{j=1}^{d_0} (-1) ^j
 \int_{S_j}
 \E\big( \det(X''_j(t) \UN_{X_{j,N}'(t)\in \ch} . \UN _{M>x}\big)
 /  X(t) =x,X'_j(t) =0\big) p_{X(t), X'_j(t)} (x,0) \sigma_j(dt)  .
\end{multline}
It is clear that for every real $x$,
$$
|p^E(x)- p_M(x)| \leq  \overline{p}(x)- p_M(x)
$$
so that the upper bounds for $ \overline{p}(x)- p_M(x)$ will automatically be upper bounds for \\
$|p^E(x)- p_M(x)|$. Moreover, as far as the authors know, no better
bounds for $|p^E(x)- p_M(x)|$ than for $ \overline{p}(x)- p_M(x)$
are known. It is an open question to determine if there exist
situations in which $p^E(x)$ is better asymptotically than
$\overline{p}(x)$.

Our next theorem gives sufficient conditions allowing to ensure that the error
$$
 \overline{p}(x)- p_M(x)
 $$
 is bounded by a Gaussian density having strictly  smaller variance than the maximum variance of the given process
 $\mathcal{ X}$ , which means that the error is super- exponentially smaller than $p_M(x)$ itself, as $ x\to +\infty$.
 In this theorem, we assume that the maximum of the variance is not attained  in
 $ S \backslash S_{d_0}$. This excludes constant variance or  some other stationary-like condition that will be
  addressed in Theorem \ref{varconst}.
As far as the authors know, the result of Theorem \ref{varpasconst}
is new even for one-parameter processes defined on a compact
interval.\medskip

 For parameter dimension $d_0>1$, the only result
of this type for non-constant variance processes of which the
authors are aware is Theorem 3.3 of \cite{TTA}.
\begin{theo} \label{varpasconst}
Assume  that the process $ \mathcal{X}$
 satisfies  conditions
 A1 -A5. With no loss of generality, we assume that $ \max_{t\in S}
 \Var(X(t)) = 1$.
In addition, we will assume that the set   $S_v$ of points $t\in
S$ where  the variance of $X(t)$ attains its maximal value is
contained in $S_{d_0} (d_0 >0)$ the non-empty face having largest
dimension and that no point  in $S_v$ is a boundary point of
$S\backslash S_{d_0}$. Then, there exist some positive constants
$C$, $\delta$ such that for every $x>0$.
\begin{equation}
  |p^E(x)- p_M(x)| \leq  \overline{p}(x)- p_M(x)\leq C \varphi( x(1+\delta)),
\end{equation}
where $\varphi(.)$ is the standard normal density.
\end{theo}

 \subsubsection*{Proof :}
Let $W$ be an open neighborhood  of the compact subset $ S_v$
 of $S$ such that $ dist( W, (S\backslash S_{d_0})) >0$ where $dist$
 denote the Euclidean distance  in $\R^d$.    For $t \in S_j \cap
 W^c$, the density
 $$
  p_{X(t), X'_j(t)} (x,0)
 $$
 can be written  as the product of the density of $X'_j(t)$ at the
 point 0, times the conditional density of $X(t)$  at the point $x$ given that $ X'_j(t) =0$,  which is Gaussian with
 some  bounded expectation and a conditional variance which is smaller  than  the unconditional
 variance, hence, bounded by some constant smaller than 1. Since
 the conditional expectations in (\ref{erdi})  are uniformly
 bounded by some constant, due to standard bounds on the moments
 of the Gaussian law, one can deduce that:
\begin{multline}\label{codix}
\overline{p}(x)- p_M(x) = \int_{W\cap S_{d_0} } \E\big(
|\det(X''_{d_0} (t))| \UN_{X_{d_0,N}'(t)\in \widehat{C}_{t,d_0}} . \UN _{M>x}
 /  X(t) =x,X'_{d_0}(t) =0\big)
 \\
  .p_{X(t), X'_{d_0}(t)} (x,0) \sigma_{d_0}(dt)
 + O( \varphi((1+\delta_1)x)),
\end{multline}
as $x \to +\infty$, for some $ \delta_1 >0$.  Our following task
is to choose $W$ such that  one can assure that the first term in
the right hand-member  of (\ref{codix}) has the same form as the
second, with a possibly different constant $\delta_1$. \\
To do this , for $s \in S $  and $t \in S_{d_0} $, let us write
the Gaussian regression formula of $X(s) $  on the pair $( X(t),
X'_{d_0} (t))$:
\begin{equation}\label{regw}
X(s) = a^t(s)  X(t) + \langle b^t(s)  , X'_{d_0}(t) \rangle + \frac{\|t-s\|^2} {2}X^t(s).
\end{equation}
where the regression coefficients $a^t(s),b^t(s)$ are respectively
real-valued and $ \R^{d_0} $-valued. \\
 From now onwards, we will only
be interested in those $t \in W$. In this case, since $W$ does
not contain boundary points  of $S\backslash S_{d_0}$, it follows
that
$$
\widehat{C}_{t,d_0} = N_{t,d_0}~\mbox{ and }
\UN_{X'_{d_0,N} (t) \in \widehat{C}_{t,d_0}} =1.
$$
Moreover, whenever $s \in S$ is close enough to $t$, necessarily, $ s
\in S_{d_0}$ and one can show that the Gaussian process $ \{
X^t(s) : t \in W \cap S_{d_0} , s\in S \}$ is bounded, in spite of
the fact  that its trajectories  are not continuous at $s=t$. For
each $t$, $ \{ X^t(s) :  s\in S \}$ is a "helix process", see
\cite{AW} for a proof of boundedness.\\
On the other hand, conditionally on $X(t) = x, X'_{d_0}(t) =0$
the event $ \{M>x\} $ can be written as
$$
\{ X^t(s) > \beta^t(s)~x, ~\mbox{ for some s}\in S \}
$$
where
\begin{equation} \label{bs1}
\beta^t(s) = \frac{2(1-a^t(s)) }{\|t-s\|^2 }.
\end{equation}

~~

 Our next goal is to prove  that if one can choose $W$ in such
 a way  that
\begin{equation}\label{szero}
\inf\{\beta^t(s) :  t \in W \cap S_{d_0}, s \in S, s \neq t\} >0,
\end{equation}
then we are done.
In fact, apply the Cauchy-Schwarz inequality  to the conditional
expectation in (\ref{codix}). Under the conditioning, the elements
of $ X''_{d_0}(t)$ are the sum of affine functions of $x$  with
bounded coefficients plus centered Gaussian variables with bounded
variances, hence, the absolute value of the conditional
expectation is bounded by an expression of the form
\begin{equation}\label{cauch}
\big( Q(t,x) \big) ^{1/2} \Big( \P\big (\sup_{s \in
S\backslash\{t\}} \frac{X^t(s) }{\beta^t(s) }>x \big) \Big)
^{1/2},
\end{equation}
 where $Q(t,x)$ is a polynomial  in $x$ of degree $2d_0$ with
 bounded coefficients. For each $t \in W\cap S_{d_0}$, the second
 factor in
 (\ref{cauch})  is bounded by
 $$
 \bigg( \P\Big( \sup\big\{\frac{X^t(s) }{\beta^t(s) } :  t \in W \cap S_{d_0}, s \in S, s \neq
 t\big\} >x \Big)\bigg)^{1/2}.
 $$
 Now, we apply to the bounded separable Gaussian process
 $$
 \Big\{\frac{X^t(s) }{\beta^t(s) } :t \in W \cap S_{d_0}, s \in S, s \neq
 t \Big\}
 $$
 the classical Landau-Shepp-Fernique inequality \cite{landau},
 \cite{fernique} which gives the bound
 $$
\P\Big( \sup\big\{\frac{X^t(s) }{\beta^t(s) } :  t \in W \cap
S_{d_0}, s \in S, s \neq
 t\big\} >x \Big) \leq C_2 \exp(- \delta_2 x^2),
 $$
 for some positive constants $ C_2, \delta_2$ and any $x>0$. Also,
 the same argument above for the density $p_{X(t), X'_{d_0}(t)} (x,0)
 $  shows  that it is bounded by a constant times the standard
 Gaussian density. To finish, it suffices to replace these bounds
 in the first term at the right-hand side of (\ref{codix}). \bigskip

 It remains to choose $W$ for (\ref{szero})  to hold true.
Consider the auxiliary process
\begin{equation}\label{defy}
  Y(s) :=\frac{ X(s)}{\sqrt{r(s,s)}},~~s\in S.
\end{equation}
Clearly, $\Var( Y(s))  =1$ for all $s \in S$. We set
$$
r^Y(s,s') := \Cov(Y(s),Y(s'))~~,~~s,s' \in S.
$$
Let us assume that $t \in S_v$. Since the function $s \Fle
\Var(X(s))$ attains its maximum value at $s=t$, it follows  that $
X(t), X'_{d_0} (t)$ are independent, on differentiation under the
expectation sign. This implies that in the regression  formula
(\ref{regw}) the coefficients are easily computed  and $ a^t(s)
=r(s,t) $ which is strictly smaller than 1 if $s \neq t $, because
of the non-degeneracy condition.

Then
\begin{equation} \label{bs2}
\beta^t(s) = \frac{2(1-r(s,t)) }{\|t-s\|^2 } \geq  \frac{2(1-r^Y(s,t)) }{\|t-s\|^2 }.
\end{equation}
Since $ r^Y(s,s) =1$ for every $s \in S$, the Taylor expansion of $r^Y(s,t)$ as a function of $s$, around
$s=t$  takes the form:
\begin{equation} \label{rho}
r^Y(s,t) = 1 + \langle s-t ,r^Y_{20,d_0}(t,t)  (s-t) \rangle + o(\|s-t\|^2),
\end{equation}
where the notation is self-explanatory. \\
Also, using that $\Var(Y(s)) =1 $ for $s \in S$, we easily obtain:
\begin{equation} \label{rho20}
-r^Y_{20,d_0,}(t,t) = \Var(Y'_{d_0} (t) ) =\Var(X'_{d_0} (t))
\end{equation}
where the last equality follows by differentiation  in (\ref{defy}) and putting $s=t$.
  (\ref{rho20}) implies that $-¨r^Y_{20,d_0,}(t,t) $ is uniformly positive definite on $ t \in S_v$, meaning that
  its minimum eigenvalue has a strictly positive lower bound. This, on account of  (\ref{bs2}) and
(\ref{rho}), already shows that
\begin{equation}\label{betapos}
\inf\{\beta^t(s) :  t \in S_v , s \in S, s \neq t\} >0,
\end{equation}
The foregoing argument also shows that
\begin{equation}\label{as}
\inf\{ -\tau (a^t)^{\prime \prime}_{d_0} (t) \tau:  t \in S_v , \tau  \in \mathcal{S}^{d_0-1}, s \neq t\} >0,
\end{equation}
since whenever $t \in S_v$, one has $ a^t(s)= r(s,t)$ so that
$$
(a^t)^{\prime \prime}_{d_0} (t) =r_{20,d_0,}(t,t).
$$
To end up, assume there is no neighborhood $W$ of $S_v$ satisfying (\ref{szero}). In that case using a compactness
argument, one can find two convergent sequences $ \{ s_n \} \subset S$ , $ \{ t_n \} \subset S_{d_0}$,
$s_n \to s_0$, $t_n \to t_0 \in S_v$ such that
$$
\beta^{t_n}(s_n)  \to \ell \leq 0.
$$
$\ell$  may be $-\infty$.\\
$t_0 \neq s_0$ is not possible, since it would imply
$$
\ell = 2 \frac{(1-a ^{t_0}(s_0))}{\|t_0-s_0\|^2} =\beta^{t_0}(s_0),
$$
which is strictly positive. \\
If $t_0 = s_0$, on differentiating in (\ref{regw}) with respect to $s$ along $ S_{d_0}$
 we get:
$$
X'_{d_0}(s) = (a^t)'_{d_0} (s) X(t)  + \langle( b^t)'_{d_0} (s)  , X'_{d_0}(t) \rangle
+  \frac{\partial_{d_0} }{\partial s}\frac{\|t-s\|^2}{2} X^t(s),
$$
where $(a^t)'_{d_0} (s)$ is a column vector of size $d_0$ and $(b^t)'_{d_0} (s)   $ is a $d_0\times d_0 $ matrix.
 Then, one must have
  $a^{t} (t) =1 $, $ (a^t)'_{d_0} (t)=0 $ . Thus
$$
\beta^{t_n}(s_n)  =  - u^T_n  (a^{t_0})^{\prime \prime}_{d_0} (t_0)  u_n + o(1),
$$
where $u_n := (s_n-t_n) /\|s_n-t_n\|$. Since $t_0 \in S_v$ we may apply (\ref{as}) and the limit $ \ell$
of $\beta^{t_n}(s_n) $ cannot be non-positive.
\hfill $\Box $\\

A straightforward application of Theorem \ref{varpasconst} is the following
\begin{cor}
Under the  hypotheses of Theorem \ref{varpasconst}, there exists positive constants $C,\delta$ such that,
for every $u>0$ :
$$
0 \leq \left| \int_u^{ +\infty} p^E(x) dx - \P(M>u)  \right| \leq \int_u^{ +\infty} \overline{p}(x) dx - \P(M>u)
\leq  C \P(\xi >u),
$$
where $ \xi$ is a centered Gaussian variable with variance $1-\delta$
\end{cor}

%

The precise order of approximation of $ \overline{p}(x) - p_M(x)$
or $ p^E(x) - p_M(x)$ as $ x \to + \infty$ remains in general an
open problem, even if one only asks for the constants $ \sigma^2_d$, $ \sigma^2_E$
respectively which govern the second order asymptotic approximation and which are defined by means of
\begin{equation}\label{logd}
\frac{1}{\sigma^2_d} := \lim_{x \to +\infty} -2 x^{ -2} \log \big[\overline{p}(x) - p_M(x)\big]
\end{equation}
and
\begin{equation}\label{logE}
\frac{1}{\sigma^2_E} := \lim_{x \to +\infty} -2 x^{ -2} \log \big|p^E(x) - p_M(x)\big|
\end{equation}
whenever these limits exist.
In general, we are unable to compute the limits (\ref{logd}) or (\ref{logE}) or even to prove that they actually exist
or differ.
 Our more general results (as well as in  \cite{AT}, \cite{TTA}) only contain lower-bounds for the liminf
 as $x \to + \infty$. This is already interesting since it gives some upper-bounds for the speed of approximation
 for $p_M(x)$ either by $\overline{p}(x)$ or $p^E(x)$.
 On the other hand, in Theorem \ref{eqisot} below, we are able to prove the existence of the limit
and compute $\sigma^2_d$ for
 a relevant class of Gaussian processes.
 \medskip

For the next theorem we need an additional condition on the
 parameter set $S$. For $S$ verifying $A1$ we define

\begin{equation} \label{capas}
 \kappa (S)=\sup_{0\leq j \leq d_0}~\sup_{t\in S_j}~\sup_{s\in S,s\neq t}\
 \frac{dist\big((t-s),C_{t,j}\big)}{\|s-t\|^2}
\end{equation}
where $dist$ is the Euclidean distance in $\R ^d$.\\
One can show that $\kappa (S)<\infty$ in each one of the following
classes of parameter sets $S$:\\
- $S$ is convex, in which case $\kappa (S)=0.$\\
- $S$ is a $C^3$ manifold, with or without boundary.\\
- $S$ verifies the following condition: For every $t \in S$ there
exists an open neighborhood $V$ of $t$ in $\R^d$ and a
$\mathcal{C}^3$ diffeomorphism $\psi :V \rightarrow B(0,r)$ (where
$B(0,r)$ denotes the open ball in $\R^d$ centered at $0$ and having
radius $r$, $r>0$) such that

$$
\psi (V \cap S)~=~C\cap B(0,r), \mbox{ where } C \mbox{ is a convex
cone}.
$$

However, $\kappa (S)<\infty$ can fail in general. A simple example
showing what is going on is the following: take an orthonormal basis
of $\R ^2$ and put
$$
S=\{(\lambda,0):0\leq \lambda \leq 1 \} \cup \{ (\mu \cos \theta
,\mu \sin \theta ):0 \leq \mu \leq 1 \}
$$
where $0<\theta < \pi$, that is, $S$ is the boundary of an angle of
size $\theta $. One easily checks that $\kappa (S)=+\infty$.
Moreover it is known \cite{AT} that in this case the EPC
approximation  does not verify  a super- exponential inequality.
More generally, sets $S$ having "whiskers" have $\kappa
(S)=+\infty$.

\medskip

\begin{theo} \label{varconst}  Let $\mathcal{X}$ be a stochastic process on $S$  satisfying
 A1 -A5.
Suppose  in addition that  $\Var(X(t) ) =1 $  for all $t \in S$ and that $\kappa (S)<+\infty$.  \\
Then
\begin{equation}\label{TTA}
\liminf _{x \to +\infty}  -2 x^{ -2} \log  \big[\overline{p}(x) - p_M(x)\big]
  \geq 1+  \inf_{t \in S} \frac{1}{ \sigma_t^2 + \overline{\lambda}(t)\kappa_t^2 }
\end{equation}
with
 $$
 \sigma_t^2 := \sup_{s\in S \backslash \{t\}} \frac{\Var\big( X(s)/X(t),X'(t)\big)}{ (1-r(s,t))^2}
 $$
 and
\begin{equation}\label{kappa}
 \kappa_t := \sup_{ s \in S \backslash \{t\}} \frac{\ dist\Big(  - \Lambda_t ^{-1} r_{01}(s,t),
 \C \Big)}{1-r(s,t)},
\end{equation}
where\begin{itemize}

  \item $\Lambda_t:=\Var(X'(t)) $
  \item $\overline{\lambda}(t)$ is the maximum eigenvalue of $ \Lambda_t$
  \item  in (\ref{kappa}), $j$ is such that $t \in S_j$ ,($j=0,1,\ldots,d_0).$
\end{itemize}

The quantity   in the right hand side of (\ref{TTA}) is strictly
bigger than $1$.
\end {theo}

\textbf{Remark}. In formula (\ref{TTA}) it may happen that the
denominator in the right-hand side is identically zero, in which
case we put $+\infty $ for the infimum. This is the case of the
one-parameter process $X(t)=\xi \cos t + \eta \sin t$ where $\xi ,
\eta$ are Gaussian standard independent random variables, and $S$ is
an interval having length strictly smaller than $\pi $.

 \subsubsection*{ Proof of Theorem \ref{varconst}}
Let us first prove that $ \sup_{t \in S } \kappa_t < \infty$.\\
For each $t \in S$, let us write the Taylor expansions
\begin{align*}
r_{01} (s,t) &= r_{01} (t,t) + r_{11} (t,t)(s-t) + O( \|s-t\|^2)
\\
 &= \Lambda_{t}(s-t)  + O( \|s-t\|^2)
\end{align*}
where $O$ is uniform on $s,t \in S$, and
\begin{align*}
1-r(s,t) &=  (s-t)^T \Lambda_{t}(s-t)  + O( \|s-t\|^2) \geq L_2
\|s-t\|^2,
\end{align*}
where $L_2$ is some positive  constant. It follows that for $s \in
S,~t \in S_j,~s \neq t$, one has:
\begin{equation}\label{courbure}
\frac{dist\Big(  - \Lambda_{t} ^{-1} r_{01}(s,t),
 \C \Big)}{1-r(s,t)}~ \leq ~L_3\frac{dist \big((t-s),\C \big)}{\|s-t\|^2}~+~L_4,
\end{equation}
where $L_3$ and $L_4$ are positive constants. So,

$$
\frac{dist\Big(  - \Lambda_{t} ^{-1} r_{01}(s,t),
 \C \Big)}{1-r(s,t)}~ \leq  ~L_3~\kappa (S)~+~L_4.
$$
which implies $ \sup_{t \in S } \kappa_t < \infty$.
\medskip\\
 With the same notations as in the proof of Theorem
\ref{varpasconst}, using (\ref{densm}) and (\ref{densmaj}), one has:

\begin{multline}\label{error6}
\overline{p}(x)- p_M(x) =   \varphi(x) \bigg[ \sum_{t \in S_0}
\E\big(\UN_{X_{t}'(t)\in \widehat{ C}_{t,0}} . \UN _{M>x}
 /  X(t) =x \big)
 \\
 + \sum_{j=1}^{d_0}
\int_{S_j} \E\big(
|\det(X^{\prime \prime}_{j} (t))| \UN_{X_{j,N}'(t)\in  \ch. \UN _{M>x}}
 /  X(t) =x,X'_{j}(t) =0\big)
 \\
(2\pi) ^{-j/2} (\det(\Var(X'_j(t))))^{-1/2}
   \sigma_{j}(dt)\bigg].
\end{multline}
Proceeding in a similar way to that of the proof of  Theorem \ref{varpasconst},
an application of the H\"older inequality  to
the conditional expectation in each term in the right-hand side of (\ref{error6})
shows that the desired result will follow as soon as we prove that:
\begin{equation}\label{liminfj}
\liminf_{x \to +\infty}  -2 x^{-2} \log  \P \big( \{X'_{j,N}  \in \ch\} \cap \{ M >x\}  /  X(t) =x,X'_j(t) =0\big)
 \geq \frac{1}{ \sigma^2_t + \overline{\lambda}(t)\kappa_t^2},
\end{equation}
for each $j= 0,1,\ldots,d_0$, where the liminf has
some uniformity in $t$.
\medskip  \\
Let us write the Gaussian regression of $X(s)$ on the pair $(X(t),X'(t))$
$$
X(s) = a^t(s) X(t) + \langle b^t(s), X'(t) \rangle + R^t(s).
$$
Since $X(t)$ and $X'(t)$ are independent, one easily computes :
\begin{align*}
a^t(s)& = r(s,t)
\\
b^t(s)& = \Lambda_t^{-1} r_{01} (s,t).
\end{align*}
Hence, conditionally on $ X(t) =x,~ X'_j(t) =0$, the events
$$
\{M >x\} ~\mbox{ and }~\{R^t(s) > (1- r(s,t)) x -r_{01}^T (s,t)
 \Lambda_t^{-1} X'_{j,N} (t) \mbox{ for some } s \in S \}
$$
coincide. \\
Denote by $(X'_{j,N} (t)|X'_j(t)=0)$ the regression of $X'_{j,N} (t)$ on $X'_j(t)=0 $.
 So, the probability in (\ref{liminfj}) can written as
\begin{equation}\label{bound1}
    \int_{\ch} \P \{\zeta^t(s) >x  -\frac{r_{01}^T (s,t)
 \Lambda_t^{-1} x'}{1-r(s,t)} \mbox{ for some } s \in S\} p _{X'_{j,N}(t)|X'_j(t)=0} (x') dx'
\end{equation}
where
\begin{itemize}
  \item $ \displaystyle \zeta ^t(s) :=\frac{ R^t(s) }{1-r(s,t) }$
  \item  $ dx'$ is the Lebesgue measure on  $ N_{t,j}$. Remember  that
  $\ch \subset N_{t,j}$.
\end{itemize}
If $ -\Lambda_t^{-1}r_{01} (s,t) \in \C$ one has
$$
-r_{01}^T (s,t)\Lambda_t^{-1} x' \geq 0
$$
for every $x'\in \ch$, because of the definition of $ \ch$.\\
If $ -\Lambda_t^{-1}r_{01} (s,t) \notin \C$, since $\C $ is a closed convex cone, we can write
$$
 -\Lambda_t^{-1}r_{01} (s,t) = z' + z^{\prime \prime}
$$
with $ z' \in \C $ , $z' \bot z^{\prime \prime}$ and $ \|z^{\prime \prime}\| =
dist( -\Lambda_t ^{-1} r_{01}(s,t), \C )$.\\
So, if $x' \in \ch$ :
$$
\frac{-r_{01}^T (s,t)\Lambda_t^{-1} x'}{1-r(s,t)} =\frac{z^{\prime T}x' +z^{\prime \prime T}x' }{1-r(s,t)}
\geq - \kappa_t \|x'\|
$$
using that $z^{\prime T}x' \geq 0 $  and the Cauchy-Schwarz inequality. It follows that in any case, if
$x' \in \ch$  the expression in (\ref{bound1}) is bounded by
\begin{equation}\label{bound2}
    \int_{\ch} \P \Big( \zeta ^t(s) >x - \kappa_t \|x'\|
  \mbox{ for some } s \in S \Big) p _{X'_{j,N} (t)|X'_j(t)=0}(x') dx'.
\end{equation}

To obtain a bound for the probability in the integrand of (\ref{bound2}) we will use the classical inequality for the tail
of the distribution of the supremum of a Gaussian process with bounded paths.

The Gaussian process $(s,t)) \Fle \zeta^t(s)$, defined on $(S \times S) \backslash \{s=t\}$ has continuous paths.
 As the pair
$(s,t)$ approches the diagonal of $S \times S$, $\zeta^t(s)$ may not
have a limit but, almost surely, it is bounded (see \cite{AW} for a
proof). (For fixed $t$, $\zeta^t(.)$ is a "helix process" with a
singularity at $s=t$, a class of processes that we have already met
above). \\

We set
\begin{itemize}
  \item $ m^t(s) := \E(\zeta^t(s)) ~~(s \neq t)$
  \item $ m := \sup_{s,t \in S ,s \neq t} | m^t(s)|$
  \item $ \mu : = \E\big( |\sup_{s,t \in S ,s \neq t}\big[\zeta^t(s) - m^t(s) \big]| \big)$.
\end{itemize}
The almost sure boundedness of the paths of $ \zeta^t(s)$ implies that $m <\infty$ and $ \mu < \infty$.
Applying the   Borell-Sudakov-Tsirelson  type inequality (see for example Adler \cite{adler90}
and references therein)
 to the centered process
$s \Fle \zeta^t(s) - m^t(s)$ defined on $ S\backslash\{t\}$
, we get whenever $x- \kappa_t \|x'\| -m - \mu >0$:
\begin{multline*}
  \P\{\zeta^t(s) > x - \kappa_t \|x'\| \mbox{ for some }s \in S\} \\
 \leq \P\{\zeta^t(s) - m^t(s) > x -\kappa_t \|x'\| -m \mbox{ for some }s \in S\}\\
  \leq  2 \exp \big( -\frac{(x- \kappa_t \|x'\| -m - \mu)^2}{2 \sigma^2_t} \big).
\end{multline*}
The Gaussian density in the integrand of (\ref{bound2}) is bounded by
$$
 (2\pi  \underline{\lambda}_j(t)) ^{\frac{j-d}{2}}
\exp \frac{\|x' - m'_{j,N}(t)\|^2}{2 \overline{\lambda}_j(t)}
$$
where $ \underline{\lambda}_j(t)$  and $\overline{\lambda}_j(t)$ are respectively  the minimum
and maximum eigenvalue of $\Var(X'_{j,N}(t)|X'_j(t))$ and $ m'_{j,N}(t)$ is the conditional expectation
$\E(X'_{j,N}(t)|X'_j(t) =0)$. Notice that $\underline{\lambda}_j(t), \overline{\lambda}_j(t), m'_{j,N}(t)$ are bounded,
$\underline{\lambda}_j(t)$ is bounded below by a positive constant and
$\overline{\lambda}_j(t) \leq \overline{\lambda}(t)$. \medskip \\
Replacing into (\ref{bound2}) we have the bound :
\begin{multline} \label{bound3}
\P \big(\{ X'_{j,N}  \in \ch\} \cap \{ M >x\}  /  X(t) =x,X'_j(t)=0\big)
\\
 \leq
 (2\pi  \underline{\lambda}_j(t))^{\frac{j-d}{2}}
 2 \int_{\ch \cap \{ x -\kappa_t \|x'\| -m - \mu >0 \}}
  \exp -\big(\frac{(x -\kappa_t \|x'\| -m - \mu)^2}{2 \sigma^2_t}
  +\frac{ \|x'-m'_{j,N} (t)\|^2}{2 \overline{\lambda}(t)}\big)~dx' \\
 + \P \Big( \| X'_{j,N}(t) |X'_j(t) =0 \| \geq \frac{x-m- \mu}{\kappa_t} \Big),
\end{multline}
where it is understood that the second term in the right-hand side
vanishes if $\kappa_t=0$. \\
Let us consider the first term in the right-hand side of (\ref{bound3}).
We have:
\begin{multline*}
\frac{(x -\kappa_t \|x'\| -m - \mu)^2}{2 \sigma^2_t}
  +\frac{ \|x'-m'_{j,N} (t)\|^2}{2 \overline{\lambda}(t)}\\
\geq
\frac{(x -\kappa_t \|x'\| -m - \mu)^2}{2 \sigma^2_t}
  +\frac{ (\|x'\|- \|m'_{j,N} (t)\|)^2}{2 \overline{\lambda}(t)}\\
= \big[ A(t) \|x'\| + B(t)(x-m-\mu) + C(t) \big] ^2  +
\frac{(x-m -\mu -\kappa_t \|m'_{j,N} (t) \| ) ^2}{2 \sigma^2_t  + 2\overline{\lambda}(t) \kappa_t^2  },
\end{multline*}
where the last inequality is obtained  after some algebra, $A(t),B(t),C(t)$ are bounded functions and
$A(t)$ is bounded below by some positive constant. \medskip
\\
So the first term in the right-hand side of  (\ref{bound3})  is bounded by :
\begin{multline} \label{bound4}
 2.(2\pi  \underline{\lambda}_j) ^{\frac{j-d}{2}}
\exp-\big( \frac{(x-m -\mu -\kappa_t \|m'_{j,N} (t) ) ^2}{2 \sigma^2_t  + 2\overline{\lambda}(t) \kappa_t^2  }\Big)
\\
\int_{R^{d-j}}  \exp -\big[\big(A(t) \|x'\| + B(t)(x-m-\mu) + C(t)\big)\big] ^2dx' \\
\leq L |x|^{d-j-1}
\exp-\Big(\frac{(x-m -\mu -\kappa_t \|m'_{j,N} (t) \|) ^2}{2 \sigma^2_t  + 2\overline{\lambda}(t) \kappa_t^2  }\Big)
\end{multline}
where $L$ is some constant. The last inequality follows easily using polar coordinates. \medskip \\
Consider now the second term in the right-hand side of  (\ref{bound3}).
Using the form of the conditional density $ p_{X'_{j,N}(t)/X'_j(t)=0}(x')$, it follows that
it is bounded by
\begin{multline} \label{bound5}
\P\big\{ \|( X'_{j,N}(t)/X'_j(t)=0 )- m'_{j,N}(t)\| \geq \frac{x- m -\mu - \kappa_t \|m'_{j,N}(t)\|}{ \kappa_t }\big \}
\\
\leq L_1 |x|^{d-j-2}
\exp-\big( \frac{(x-m -\mu -\kappa_t \|m'_{j,N} (t) \|) ^2}{  2\overline{\lambda}(t) \kappa_t^2  }\Big)
\end{multline}
where $L_1$ is some constant. Putting together (\ref{bound4}) and
(\ref{bound5}) with (\ref{bound3}), we obtain (\ref{liminfj}).
\hfill $\Box $\newline

The following two corollaries  are straightforward consequences of Theorem \ref{varconst}:
\begin{cor}
Under the hypotheses of Theorem \ref{varconst} one has
$$
\liminf_{x \to +\infty} -2 x^{-2} \log | p^E(x) -p_M(x) | \geq 1 +
\inf_{t \in S} \frac{1}{\sigma^2_t  + \overline{\lambda}(t)
\kappa_t^2 }.
$$
\end{cor}
\begin{cor} Let $\mathcal{X}$ a stochastic process on $S$  satisfying
 A1 -A5.
Suppose  in addition that $\E(X(t)) =0$,  $\E(X^2(t) ) =1 $, $\Var(X'(t) = I_d$   for all $t \in S$.\\
Then
$$
\liminf_{u \to +\infty} -~2 u^{-2} \log \Big| \P(M>u)  -
\int_u^{+\infty} p^E(x) dx \Big| \geq  1 + \inf_{t \in S}
\frac{1}{\sigma^2_t  + \kappa_t^2 }.
$$
and
$$
p^E(x) = \Big[ \sum_{j=0} ^{d_0} (-1)^j  (2 \pi) ^{-j/2} g_j\overline{ H}_j(x) \Big] \varphi(x).
$$
where $ g_j$ is given by (\ref{geom}) and $\overline{ H}_j(x)$ has
been defined in Section \ref{comppbar}.
\end{cor}

The proof follows directly from Theorem \ref{varconst} the
definition of $p^E(x)$ and the results in \cite{A81}.

\section{Examples}

 1)
A simple application  of Theorem \ref{varpasconst} is the following.
Let $\mathcal{X}$ be a one parameter real-valued
centered
Gaussian process with regular paths, defined on the interval $[0,T]$ and
 satisfying an adequate non-degeneracy condition.
Assume that the variance $v(t)$ has a unique maximum, say
$1$ at the interior point $t_0$, and $k=\min\{j:v^{(2j)}(t_0)\neq 0\}<\infty$.
Notice that $v^{(2k)}(t_0)<0$.
Then, one can obtain the equivalent of $p_M(x)$ as $x\to\infty$ which is given by:
\begin{equation}\label{equiv}
p_M(x)\simeq\frac{1-v^{\prime \prime}(t_0)/2}{kC_k^{1/k}}
                                          \E\left(|\xi|^{
                                                         \frac{1}{2k}          -1}
                                                         \right) x^{1-1/k}\varphi(x),
\end{equation}
where $\xi$ is a standard normal random variable and
$C_k=-\frac{1}{(2k)!}v^{(2k)}(t_0)+\frac{1}{4}[v^{\prime
\prime}(t_0)]^2 \UN_{k=2}$. The proof is a direct application of the
Laplace method. The result is new for the density of the maximum,
but if
 we integrate the density from $u$ to $+\infty$, the corresponding  bound for $ \P\{M >u\}$ is known under weaker
  hypotheses (Piterbarg \cite{piter96}).
 \medskip

 2) Let the process $\mathcal{ X}$ be centered and satisfy A1-A5. Assume that the the law
 of the process is isotropic
and stationary, so that the covariance has the form (\ref{isotro})
and verifies the regularity condition of Section \ref{comppbar}. We
add the simple normalization $ \rho' = \rho'(0) = -1/2$. One can
easily check that
\begin{equation}\label{sigisot}
\sigma^2 _t =\sup_{s\in S \backslash\{t\}} \frac{1-\rho^2(\|s-t\|^2) -4 \rho^{\prime 2}(\|s-t\|^2)\|s-t\|^2}
{ [1-\rho(\|s-t\|^2)]^2}
\end{equation}
Furthermore if
\begin{equation}\label{ropneg}
  \rho'(x) \leq 0 \mbox{ for } x \geq 0
\end{equation}
one can show that the sup in (\ref{sigisot}) is attained as $\|s-t\| \to 0$ and is independent of $t$. Its value is
$$
\sigma^2 _t = 12 \rho^{\prime \prime}-1.
$$
The proof is elementary (see \cite{ABW} or \cite{TTA}).\\

Let $S$ be a convex set. For t $ \in S_j$, $s\in S$:
\begin{equation}\label{distance}
dist\big(-r_{01} (s,t) ,\C\big) = dist\big( -2\rho'(\|s-t\|^2) (t-s) ,\C\big).
\end{equation}
The convexity of $S$ implies that $(t-s) \in \C$. Since $\C$ is a
convex cone and $-2\rho'(\|s-t\|^2) \geq 0$, one can conclude  that
$-r_{01} (s,t) \in \C$ so that the distance in (\ref{distance}) is
equal to zero. Hence,
$$
\kappa_t =0 \mbox{ for every } t \in S
$$
and an application of  Theorem \ref{varconst} gives the inequality
\begin{equation}\label{isotcon}
  \liminf_{x \to +\infty} -\frac{2}{x^2}  \log\big[ \overline{p}(x)-p_M(x)  \big]  \geq 1 +\frac{ 1}{12 \rho^{\prime \prime} -1}.
\end{equation}

A direct consequence is that the same inequality  holds true when replacing $\overline{p}(x)-p_M(x)$ by
$ |p^E(x)-p_M(x)|$ in (\ref{isotcon}), thus obtainig the main explicit example in
Adler and Taylor \cite{AT}, or  in Taylor et al. \cite {TTA}.

Next, we improve (\ref{isotcon}). In fact, under the same hypotheses, we prove that the liminf is an ordinary limit
 and the sign $\geq$ is an equality sign. We state this as
 \begin{theo} \label{eqisot}
Assume that $\mathcal{X}$ is centered, satisfies hypotheses A1-A5,
the covariance has the form (\ref{isotro}) with $\rho' (0)=-1/2,~
\rho' (x) \leq 0~for~~x \geq 0$. Let $S$ be a convex set, and $d_0 =
d \geq 1$. Then
\begin{equation}\label{eqisote}
\lim_{x \to +\infty} -\frac{2}{x^2}  \log\big[ \overline{p}(x)-p_M(x)  \big]  = 1 +\frac{ 1}{12 \rho^{\prime \prime} -1}.
\end{equation}
\end{theo}
{\bf Remark} Notice that since $S$ is convex,  the added
hypothesis that the maximum dimension $d_0$ such that $S_j$ is not
empty is equal to $d$ is not an actual restriction.
\subsubsection*{Proof of Theorem \ref{eqisot}}

In view of (\ref{isotcon}), it suffices to prove that
\begin{equation}\label{isotcon2}
\limsup_{x \to +\infty} -\frac{2}{x^2}  \log\big[
\overline{p}(x)-p_M(x)  \big]  \leq 1 +\frac{ 1}{12 \rho^{\prime \prime} -1}.
\end{equation}
Using (\ref{densm}) and the definition of $ \overline{p}(x) $ given
by (\ref{densmaj}), one has the inequality
\begin{equation}\label{pas1}
\overline{p}(x)-p_M(x) \geq  (2 \pi)^{ -d/2} \varphi(x)
\int_{S_d} \E\big( |\det(X^{\prime \prime}(t) ) |\UN _{M>x} /X(t) =x, X'(t) =0) \sigma_d(dt),
\end{equation}
where our lower bound only contains the term corresponding to the largest dimension and we have already
 replaced  the density $p_{X(t),X'(t)} (x,0) $ by its explicit expression  using the law of the process.
 Under the condition  $\{ X(t) =x, X'(t) =0\} $ if $ v_0^T X^{\prime \prime}(t) v_0 >0 $ for some $v_0 \in \mathcal{S}^{d-1}$,
  a Taylor expansion implies that $M>x$. It follows that
\begin{multline}\label{pas2}
\E\big( |\det(X^{\prime \prime}(t) )| \UN _{M>x} /X(t) =x, X'(t) =0\big) \\
\geq
\E\big( |\det(X^{\prime \prime}(t) ) |\UN _{ \displaystyle \sup_{v \in \mathcal{S}^{d-1}} v^T X^{\prime \prime}(t) v >0} /X(t) =x ,X'(t) =0\big).
\end{multline}
We now apply Lemma \ref{formules} which describes the conditional distribution of $X^{\prime \prime}(t)$  given
$X(t) =x, X'(t) =0$ . Using the  notations  of this lemma, we may write the right-hand side of (\ref{pas2}) as :
$$
\E\big(| \det(Z- x Id )| \UN _{ \displaystyle \sup_{v \in
\mathcal{S}^{d-1}} v^T Z v >x}\big),
$$
which is obviously bounded below by
\begin{multline} \label{pas3}
\E\big(| \det(Z- x Id )| \UN _{Z_{11} >x}\big)\\
= \int_x^{+\infty} \E\big(| \det(Z- x Id )| /Z_{11} =y \big)  (2 \pi) ^{-1/2} \sigma ^{-1}
\exp\big( - \frac{y^2}{2 \sigma^2}\big) dy,
\end{multline}
where $\sigma^2 := \Var(Z_{11} ) = 12 \rho^{\prime \prime}-1$. The conditional distribution of $Z$ given $Z_{11} =y$ is easily
deduced from Lemma \ref{formules}. It can be represented by the random $d \times d$ real symmetric matrix
$$
\widetilde{Z} := \left(
\begin{array}{ccccc}
y& Z_{12}&\ldots &\ldots&Z_{1d}     \\
 & \xi_2 + \alpha y &Z_{23} & \ldots&Z_{2d}     \\
  & & \ddots & &\\
  &&&&\xi_d + \alpha y
  \end{array}
  \right) ,
 $$
 where the random variables $\{\xi_2, \ldots, \xi_d,Z_{ik}, 1 \leq i<k \leq d \}$
 are independent centered Gaussian with
 $$
 \Var(Z_{ik}) = 4 \rho^{\prime \prime} ~~(1 \leq i<k \leq d)~~ ;~~ \Var(\xi_i) = \frac{16 \rho^{\prime \prime} (8 \rho^{\prime \prime}-1)}{12 \rho^{\prime \prime} -1}
~~(i=2,\ldots,d)~~;~~ \alpha =\frac{4 \rho^{\prime \prime}-1}{12 \rho^{\prime \prime} -1}
$$
Observe that $0 <\alpha <1$.\\
Choose now $\alpha_0$ such that $(1+\alpha_0) \alpha <1$. The expansion of $ \det( \widetilde{Z}-x Id)$ shows that
if $  x(1+\alpha_0) \leq y \leq x(1+\alpha_0) +1$ and $x$ is large enough, then
$$
\E\big(| \det(\widetilde{Z}- x Id )| \big) \geq L ~\alpha_0 ( 1-\alpha (1+\alpha_0))^{d-1}~x^d,
$$
where $L$ is some positive constant. This implies that
$$
 \frac{1}{\sqrt{2 \pi} \sigma} \int_x^{+\infty} \exp(- \frac{y^2}{2\sigma^2})\E\big(| \det(\widetilde{Z}- x Id )| \big) dy
 \geq
\frac{L}{\sqrt{2 \pi} \sigma} \int_{x(1+\alpha_0)}^{x(1+\alpha_0) +1} \exp(- \frac{y^2}{2\sigma^2})
\alpha_0 ( 1-\alpha (1+\alpha_0))^{d-1}~x^d dy
$$
for $x$ large enough. On account of
(\ref{pas1}),(\ref{pas2}),(\ref{pas3}), we conclude that for $x$
large enough,
$$
\overline{p}(x)-p_M(x) \geq  L_1 x^d  \exp -\Big[ \frac{x^2}{2} +\frac{ ( x(1+\alpha_0) +1)^2}{2 \sigma^2}\Big].
$$
for some new positive constant $L_1$. Since $\alpha_0$ can be chosen
arbitrarily small,
 this implies (\ref{isotcon2}). \\
 \hfill $\Box $ \medskip

3) Consider the same  processes  of  Example 2, but now defined on
the non-convex set $ \{a \leq \|t\| \leq b\}$,
 $0<a<b$. The same calculations as above  show that $ \kappa_t=0$ if $a<\|t\| \leq b$  and
$$
\kappa_t =  \max \Big\{  \sup_{z \in[2a, a+b]} \frac{ -2 \rho'(z^2) z}{ 1-\rho(z^2) },
\sup_{\theta \in[0,\pi]}\frac{ -2a \rho'( 2a^2( 1-\cos \theta ) )( 1-cos \theta )}{ 1-\rho(2a^2( 1-\cos \theta )) }
\Big\},
$$
for $\|t\| =a$.\medskip

    4) Let us keep the same hypotheses as in Example 2 but without assuming that the covariance
     is decreasing as in (\ref{ropneg}).
The variance is still given by (\ref{sigisot}) but $ \kappa_t$ is not necessarily equal to zero.
    More precisely, relation (\ref{distance}) shows that
    $$
    \kappa_t \leq  \sup_{s \in S\backslash \{t\} } 2 \frac{\rho'(\| s-t \| ^2 )^+ \|s-t\|}{1-\rho(\|s-t\|^2)}
    $$
    The normalization: $ \rho' = -1/2$ implies that the process $\mathcal{ X}$ is "identity speed",
     that is \\
      $\Var(X'(t))= I_d$
 so that $\overline{\lambda}(t) =1$. An application of  Theorem \ref{varconst}  gives
 \begin{equation}\label{xxxx}
  \liminf_{x \to +\infty} -\frac{2}{x^2}  \log\big[ \overline{p}(x)-p_M(x)  \big]  \geq
  1 + 1/Z_\Delta.
\end{equation}
where
$$
Z_\Delta := \sup_{z \in (0,\Delta]} \frac{1-\rho^2(z^2) -4 \rho^{\prime 2}(z^2)z^2}
{ [1-\rho(z^2)]^2} + \max_{z \in (0,\Delta]} \frac{4 \big[\rho'(z ^2 )^+ z\big]^2}{[1-\rho(z^2)]^2},
$$
and $\Delta$ is the diameter of $S$.\\

 5)   Suppose that

\begin{itemize}
  \item the process $\mathcal{X}$ is stationary with
covariance $ \Gamma(t) := \Cov ( X(s),X(s+t))$  that satisfies
  $\Gamma(s_1,\ldots,s_d)=  \prod _{i=1,...,d} \Gamma_i(s_i)$
  where $\Gamma_1,..., \Gamma_d$ are $d$ covariance functions on $\R$ which are
  monotone, positive on
  $[0, + \infty)$ and of class
  $ \mathcal{C}^4$,
  \item $S$ is a rectangle
  $$
  S = \prod _{i=1,...,d} [a_i,b_i] ~~,a_i<b_i.
  $$
  \end{itemize}
  Then, adding an appropriate non-degeneracy condition, conditions A2-A5
  are fulfilled and Theorem  \ref{varconst} applies

It is easy to see that
  $$
 -  r_{0,1} (s,t) = \left[ \begin{array}{c}
 \Gamma'_1(s_1-t_1)  \Gamma_2(s_2-t_2) \ldots  \Gamma_d(s_d-t_d) \\
 \vdots  \\
\Gamma_1(s_1-t_1) \ldots \Gamma_{d-1}(s_{d-1}-t_{d-1}). \Gamma'_d(s_d-t_d)\\
\end{array} \right]
$$
belongs to $\C$ for every $s \in S$. As a consequence $\kappa_t=0$
for all $t \in S$. On the other hand, standard regressions formulae
show that
$$
\frac{\Var\big( X(s)/X(t),X'(t)\big)}{ (1-r(s,t))^2} =
\frac{1 -\Gamma_1^2\ldots \Gamma_d^2-
\Gamma^{\prime 2}_1 \Gamma_2^2\ldots \Gamma_d^2 - \cdots-\Gamma^{ 2}_1 \ldots \Gamma_{d-1}^2
\Gamma_d^{\prime 2  }}{ (1-\Gamma_1\ldots\Gamma_d)^2},
$$
where $\Gamma_i$ stands for $\Gamma_i(s_i-t_i)$. Computation  and
maximisation of $\sigma^2_t$ should be performed numerically in each
particular case.

\end{document}